\documentclass[12pt]{amsart}
\usepackage{amscd,amssymb,epsfig}
\usepackage{graphics}

\newtheorem{theorem}{Theorem}[section]
\newtheorem{lemma}[theorem]{Lemma}

\newtheorem{corollary}[theorem]{Corollary}
\theoremstyle{definition}

\theoremstyle{remark}

\numberwithin{figure}{section}
\numberwithin{table}{section}

\makeatletter

\@addtoreset{equation}{section}
\makeatother

\newcommand{\Hom}{\operatorname{Hom}}

\newcommand\bC{{\mathbb C}}

\newcommand\bR{{\mathbb R}}
\newcommand\bZ{{\mathbb Z}}
\newcommand\Cg{{{\mathbb C}_g}}

\newcommand\dc{\overline{\partial}}
\newcommand\dd{{\partial}}

\newcommand{\dvol}{\operatorname{dvol}}
\newcommand\fiber{{\text{fiber}}}
\newcommand\harmonic{{\mathcal{H}}}
\newcommand\hatPhi{{\widehat{\Phi}}}
\newcommand\Hg{{\mathfrak{H}_g}}
\newcommand\Hodge{{\Lambda_{\mathbb{M}_g}}}
\newcommand\hyp{\text{hyperbolic}}

\newcommand\Jac{{\operatorname{Jac}}}

\newcommand\Mg{{\mathcal{M}_g}}
\newcommand\Mgone{{\mathcal{M}_{g, 1}}}
\newcommand\Mgstar{{\mathcal{M}_{g, *}}}
\newcommand\moduli{{{\mathbb M}_g}}
\newcommand\mudot{{\stackrel{\centerdot}{\mu}}}
\newcommand\omegaone{{\omega_{(1)}}}
\newcommand\omegaonedot{{\stackrel{\centerdot}{\omega_{(1)}}}}

\newcommand{\Res}{\operatorname{Res}}

\newcommand\stardot{{\stackrel{\centerdot}{\ast}}}
\newcommand{\Sym}{\operatorname{Sym}}
\newcommand{\Symp}{\operatorname{Sp}}

\newcommand\TCM{{T^\times_{\mathbb{C}_g/\mathbb{M}_g}}}
\newcommand\Tg{{\mathcal{T}_g}}
\newcommand{\trace}{\operatorname{trace}}
\newcommand\WP{\mbox{\tiny WP}}

\begin{document}

\title{Canonical 2-forms on the moduli space of Riemann surfaces}

\author{Nariya Kawazumi}
\thanks{Work partially supported by 
Grant-in-Aid for Scientific Research (A) (No.18204002), 
the Japan Society for Promotion of Sciences.
}
\address{
Department of Mathematical Sciences,\\
University of Tokyo,\\
Tokyo, 153-8914 Japan\\
email:\,\tt{kawazumi@ms.u-tokyo.ac.jp}
}

\begin{abstract} 
As was shown by Harer \cite{H1} \cite{H2}, the second homology of ${\mathbb M}_g$, the moduli space of compact Riemann surfaces of genus $g$, is of rank $1$, provided $g \geq 3$. This means there exists 
a nontrivial second de Rham cohomology class on ${\mathbb M}_g$ which is unique 
up to a constant factor. But several canonical $2$-forms 
on the moduli space have been constructed in various geometric contexts, and they differ from each other. 
In this article we review 
some constructions of such canonical $2$-forms 
in order to provide material 
for future research on the ``secondary geometry" 
of the moduli space ${\mathbb M}_g$. 
\end{abstract}



\maketitle

\tableofcontents

\section{Introduction}\label{s-1}

Let $g \geq 2$ be an integer. 
The moduli space of compact Riemann surfaces of 
genus $g$, $\moduli$, is the quotient space of  Teichm\"uller 
space $\Tg$ by the natural action of the mapping class group $\Mg$,
$\moduli = \Tg/\Mg$. Since Teichm\"uller space is contractible, 
the real cohomology of the mapping class group is 
isomorphic to that of the moduli space. 
As was shown by Harer \cite{H1} \cite{H2}, the second homology of ${\mathbb M}_g$ is of rank $1$ if $g \geq 3$. This means there exists 
a nontrivial second de Rham cohomology class on ${\mathbb M}_g$ which is unique 
up to a constant factor. But several canonical $2$-forms 
on the moduli space have been constructed in various geometric contexts, and they differ from each other. 
In this article we review 
some constructions of such canonical $2$-forms 
in order to provide material 
for future research on the ``secondary geometry" 
of the moduli space ${\mathbb M}_g$. 
\par

The signature of the total space of a fiber bundle 
is not necessarily equal to the product of the 
signatures of the base space and the fiber. 
The first example for this phenomenon was
given by Kodaira \cite{Ko} and Atiyah \cite{At},  who constructed a certain branched covering space
of the product of two compact Riemann surfaces. 
The covering space has non-zero signature, while 
the signature of any compact Riemann surface is zero. 
We may regard the covering space as a family of compact Riemann 
surfaces parametrized by a compact Riemann surface, so that it defines 
a non-trivial $2$-cycle on the space $\moduli$. 
As was formulated by Meyer 
\cite{Me1} \cite{Me2}, the signature of the total space 
of a family of compact Riemann surfaces defines 
a non-trivial $2$-cocycle of the mapping class group 
$\Mg$ and this provides a non-trivial cohomology class of degree $2$ on the space $\moduli$. 
Nowadays this cocycle is called the Meyer cocycle 
and it has been playing an essential role in the topological study of fibered complex surfaces. 
See \cite{AE} and \cite{AKonno} for details. \par
The first and the second Betti numbers of the space 
$\moduli$, or equivalently, those of the group $\Mg$,
are given by
\begin{eqnarray}
& b_1(\moduli) = 0, & \mbox{\cite{Mu1} \cite{Po} \cite[p.223]{H1}}\label{b1}\\
& b_2(\moduli) = 1, \quad \mbox{if $g \geq 3$.} 
& \mbox{\cite{H1} \cite{H2}}\label{b2}
\end{eqnarray}
For alternative computations of $b_2(\moduli)$, see 
\cite{AC} \cite{KS} \cite{Pi}. The group
$H^2(\moduli; \bR)$ is generated by the cohomology class of the Meyer cocycle. 
In the case $g=2$ we have $b_2({\mathbb M}_2) = 0$ because of Igusa's result ${\mathbb M}_2 = \bC^3/(\bZ/5) \simeq *$ \cite{JI}.
\par
Mumford \cite{Mu} and Morita \cite{Mo1} independently introduced a series of cohomology classes $e_n = (-1)^{n+1}\kappa_n \in H^{2n}(\moduli)$, $n \geq 1$, 
the Morita-Mumford classes or the tautological classes.
They are defined as follows. Let $\pi: \Cg \to \moduli$ 
be the universal family of compact Riemann surfaces 
of genus $g$. The relative tangent bundle of the map $\pi$,
$T_{\mathbb{C}_g/\mathbb{M}_g}$, the kernel of the differential $d\pi:
T\Cg \to 
\pi^*T\moduli$, is a complex line V-bundle over $\Cg$. 
The $n$-th Morita-Mumford class 
$e_n = (-1)^{n+1}\kappa_n$, $n\geq 1$, is defined to be the 
integral of the $(n+1)$-st power of the Chern class of the bundle 
$\TCM$ along the fiber
\begin{equation}
e_n = (-1)^{n+1}\kappa_n = \int_\fiber c_1(\TCM)^{n+1} \in 
H^{2n}(\moduli).
\end{equation}
The first one $e_1 = \kappa_1$ is $3$ times the cohomology class of the Meyer cocycle. 
As was proved by Morita \cite{Mo2} and Miller \cite{Mi}, 
the Morita-Mumford classes are algebraically independent in the stable range $* < \frac23g$ \cite{H3} of the cohomology algebra $H^*(\moduli; \bR)$. Their proofs generalize the construction of Kodaira and Atiyah.
In 2002 Madsen and Weiss \cite{MW} proved that the cohomology algebra $H^*(\moduli; \bR)$ in the 
stable range is generated by the Morita-Mumford classes. \par
From the results (\ref{b1}) and (\ref{b2}) the simplest non-trivial cohomology classes on $\moduli$ are of degree $2$, and they are 
unique up to a constant factor. 
But several $2$-forms on $\moduli$, or equivalently $\Mg$-equivariant $2$-forms on Teichm\"uller space $\Tg$, have been canonically constructed in various geometric contexts. 
\par
From the uniformization theorem any compact Riemann surface $C$ of genus $g \geq 2$ admits a unique hyperbolic 
metric. The volume form of the hyperbolic metric defines 
the Weil-Petersson pairing on the cotangent space 
$T^*_{[C]}\moduli$ involved with no addditional information. 
As was shown by Wolpert \cite{Wo} the Weil-Petersson-K\"ahler form $\omega_{\WP}$ represents the first 
Morita-Mumford class $e_1$. Thus we obtain a canonical 
$2$-form representing $e_1$.\par
The period map is a canonical map defined on Teichm\"uller space into the Siegel upper halfspace $\Hg$. We have a canonical $2$-form on $\Hg$ whose pullback represents 
the class $e_1$ on the moduli space $\moduli$. \par
We have another canonical metric on a compact Riemann surface. A natural Hermitian product on the space of holomorphic $1$-forms defines the volume form $B$ in 
\ref{B} which induces a Hermitian metric on the Riemann surface. The Arakelov-Green function is derived from the volume form $B$. As will be stated in \S\ref{s-hme} and 
\S\ref{s-fiber}, a higher analogue of the period map is 
constructed and yields other canonical $2$-forms representing $e_1$. These forms are closely related to the volume form $B$.\par
All of them differ from each other. 
As to $2$-forms representing non-trivial cohomology classes of degree $2$ on the moduli space $\moduli$, the term `canonical' does not imply `unique'. 
The difference of such forms should induce some secondary object on the moduli space $\moduli$. 
Assume $g \geq 3$. If we have two real $(1, 1)$-forms 
$\psi_1$ and $\psi_2$ on $\moduli$ representing $e_1$, 
then there exists a real-valued function $f \in C^\infty(M; \bR)$ such that $
\psi_2 - \psi_1 = \frac{\sqrt{-1}}{2\pi}\dd\dc f$. 
Such a function $f$ is unique up to a constant.
See Lemma \ref{d-d}. This function captures the difference between these two forms, so that it should describe a certain relation between the two geometric contexts behind these forms. 
\par

In this article we review some constructions of canonical $2$-forms. 
In \S\ref{s-cotang} we give a short review on the cotangent spaces of moduli spaces. They are naturally isomorphic to some spaces of quadratic differentials. In \S\ref{s-wp} we take a quick glance at the Weil-Petersson K\"ahler form, which is related to the Virasoro cocycle through the Krichever construction. 
The most classical $2$-form on $\moduli$ is the pullback of the first Chern form on the Siegel upper halfspace $\Hg$ by the period map $\Jac$, or equivalently the first Chern form of the Hodge bundle on $\moduli$. We explain this form in \S\S\ref{s-chern} and \ref{s-rau}. The Hodge bundle yields all the odd Morita-Mumford classes but not the even ones.
We can obtain other canonical differential forms on the moduli space representing all the Morita-Mumford class $e_i$, $i \geq 1$, through a higher analogue of the period map, 
and this is described in \S\S \ref{s-earle} and \ref{s-hme}. 
Among them some $2$-forms seem to be related to  Arakelov geometry, as
will be discussed in \S \ref{s-fiber}. \par 
First of all the author
thanks Athanase Papadopoulos for careful reading and valuable comments on
this article. Furthermore he thanks Leon Takhtajan for helpful comments
on an earlier version. 
\par

\section{The cotangent space of the moduli space}\label{s-cotang}

Let $C$ be a compact Riemann surface of genus $g \geq 2$, $P_0$ a point on $C$. Then
we denote by $H^q(C; aK + bP_0)$, $q = 0, 1$, and 
$a, b \in \bZ$, 
the $q$-th cohomology group $H^q(C; \mathcal{O}_C(
T^*C^{\otimes a}\otimes [P_0]^{\otimes b}))$. 
Moreover we denote by $\Omega^q(C)$ the complex-valued
$q$-currents on $C$ for $0 \leq q \leq 2$ . 
The Hodge $*$-operator
 $*: \left(T^*_\bR C\right)\otimes\bC 
\to \left(T^*_\bR C\right)\otimes \bC$ on the cotangent 
bundle of $C$ depends only on the complex structure of $C$.
The $-\sqrt{-1}$-eigenspace is the holomorphic cotangent bundle
$T^*C$, and the $\sqrt{-1}$-eigenspace is the antiholomorphic
cotangent bundle $\overline{T^*C}$. 
The operator $*$ decomposes the space $\Omega^1(C)$ into 
the $\pm \sqrt{-1}$-eigenspaces
$$
\Omega^1(C) = \Omega^{1, 0}(C) \oplus \Omega^{0, 1}(C),
$$ 
where $\Omega^{1, 0}(C)$ is the $-\sqrt{-1}$-eigenspace and 
$\Omega^{0, 1}(C)$ the $\sqrt{-1}$-eigenspace. 
Throughout this article we denote by $\varphi'$ and $\varphi''$ 
the $(1, 0)$- and the $(0, 1)$-parts of $\varphi \in \Omega^1(C)$, 
respectively, i.e., 
$$
\varphi = \varphi' + \varphi'', \quad 
*\varphi = -\sqrt{-1}\varphi' + \sqrt{-1}\varphi''.
$$
If $\varphi$ is harmonic, then $\varphi'$ is holomorphic and 
$\varphi''$ anti-holomorphic.\par

The Kodaira-Spencer map gives a natural isomorphism
\begin{equation}
T_{[C]}\moduli = H^1(C; -K).\label{KS}
\end{equation}
To look at the isomorphism (\ref{KS}) more explicitly, consider 
a $C^\infty$ family of compact Riemann surfaces
$C_t$, $t \in \bR$, $\vert t \vert \ll 1$, with 
$C_0 = C$. 
The family $\{C_t\}$ is trivial as a $C^\infty$ 
fiber bundle over an interval near $t = 0$, 
so that we have a $C^\infty$ family of $C^\infty$ diffeomorphisms 
$f^t: C \to C_t$
with $f^0 = 1_C$. 
In general, if $\bigcirc = \bigcirc_t$ is a ``function" in 
$t \in \bR$, $\vert t \vert \ll 1$, then we write simply 
\begin{equation*}
\mathop\bigcirc^{\centerdot} = 
\frac{d}{dt}\Bigr\vert_{t=0}\bigcirc_t.
\end{equation*}
For example, we denote 
$$
\mudot = \frac{d}{dt}\Bigr\vert_{t=0}\mu(f^t).
$$
Here $\mu(f^t)$ is the complex dilatation of the diffeomorphism $f^t$. 
Let $z_1$ be a complex coordinate on $C$, and $\zeta_1$ on $C_t$. 
The complex dilatation
 $\mu(f^t)$ is defined locally by 
\begin{equation*}
\mu(f^t) = \mu(f^t)(z_1)\frac{d}{dz_1}\otimes d\overline{z_1} = 
\frac{(\zeta_1\circ
f^t)_{\overline{z_1}}}{(\zeta_1\circ
f^t)_{z_1}}\frac{d}{dz_1}\otimes d\overline{z_1},
\end{equation*}
which does not depend on the choice of the coordinates $z_1$ and $\zeta_1$. 
The Dolbeault cohomology class $[\mudot] \in H^1(C; -K)$ is exactly 
the tangent vector $\frac{d}{dt}\bigr\vert_{t=0}[C_t] \in T_{[C]}\moduli$. \par
We define a linear operator $S = S[\mudot]: \Omega^1(C) \to \Omega^1(C)$ by 
$$
S(\varphi) = S(\varphi') + S(\varphi'') :=
-2\varphi'\mudot - 2\varphi''\overline{\mudot},
$$
for $\varphi = \varphi' + \varphi''$,  $\varphi' \in \Omega^{1, 0}(C)$, 
$\varphi'' \in \Omega^{0, 1}(C)$. 
From straightforward computation we have 
\begin{equation}
\stardot = *S = -S*: \Omega^1(C) \to \Omega^1(C).
\end{equation}

By Serre duality we have a natural isomorphism
\begin{equation}
T^*_{[C]}\moduli = H^0(C; 2K).
\end{equation}
The space $H^0(C; 2K)$ consists of the holomorphic quadratic 
differentials on $C$. For any holomorphic quadratic 
differential $q$ the covariant tensor $q\mudot$ can be regarded 
as a $(1, 1)$-form on $C$. The integral 
$\int_Cq\mudot$ is just the value of the covector $q$ at the 
tangent vector $[\mudot] = \frac{d}{dt}\bigr\vert_{t=0}[C_t]$. \par
Let $\Cg$ denote the moduli space of pointed compact Riemann 
surfaces $(C, P_0)$ of genus $g$ with $P_0 \in C$. The forgetful map 
$\pi: \Cg \to \moduli$, $[C, P_0] \mapsto [C]$, 
can be interpreted as the universal family of 
compact Riemann surfaces on the moduli space $\moduli$. 
We identify 
\begin{equation}
T_{[C, P_0]}\Cg = H^1(C; -K-P_0), \quad \mbox{and}\quad
T^*_{[C, P_0]}\Cg = H^0(C; 2K+P_0)
\end{equation}
in a way similar to the space 
$\moduli$.\par
The relative tangent bundle of the forgeful map $\pi$
with the zero section deleted
$$
\TCM = T_{\mathbb{C}_g/\mathbb{M}_g} \setminus
\mbox{(zero section)}
$$
can be interpreted as the moduli space of 
triples $(C, P_0, v)$ of  genus $g$. Here 
$C$ is a compact Riemann surface of genus $g$, 
$P_0 \in C$, and 
$v \in T_{P_0}C\setminus\{0\}$. 
Similarly the space of quadratic differentials $H^0(C; 2K+2P_0)$ 
is identified with the cotangent space of $\TCM$
\begin{equation}
T^*_{[C, P_0, v]}\TCM = H^0(C; 2K+2P_0).
\end{equation}
Moreover this space is closely related to Ehresmann connections on
the bundle $T_{\mathbb{C}_g/\mathbb{M}_g}$.
In general, let $\varpi: L \to M$ be a holomorphic line bundle 
over a complex manifold $M$, and $L^\times$ the total space with 
the zero section deleted $L^\times = L \setminus \mbox{(zero section)}$. 
We denote by $R_a$ the right action of $a \in \bC^\times := \bC \setminus \{0\}$ on the space $L^\times$, and by $Z$ the vector field
on $L^\times$ generated by the action $R_a$
$$
Z := \frac{d}{dt}\Bigr\vert_{t=0}R_{e^t}.
$$
An Ehresmann connection $A$ (of type $(1, 0)$)
 on the bundle $L$ 
is a $(1, 0)$-form on the space $L^\times$ with the conditions
\begin{eqnarray*}
& A(Z) = 1, \quad\mbox{and}\\
& {R_{e^t}}^*A = A, \quad \forall t \in \bR
\end{eqnarray*}
\cite{Atconn}\cite{KN}. 
In other words, it is a splitting of the extension of holomorphic vector bundles over $M$
\begin{equation*}
0 \to T^*M \mathop{\to}^{\varpi^*} (T^*L^\times)/\bC^\times 
\mathop{\to}^{Z} \bC \to 0.
\end{equation*}
Then there exists a unique $(1, 1)$-form $c_1(A)$ 
on $M$ such that $\frac{\sqrt{-1}}{2\pi}dA = \varpi^*c_1(A)$. 
The form $c_1(A)$ is, by definition, the Chern form of the connection 
$A$ and represents the first Chern class of the line bundle $L$
$$
[c_1(A)] = c_1(L) \in H^2(M; \bR).
$$
\par
Now we let $M = \Cg$ and $L = T_{\mathbb{C}_g/\mathbb{M}_g}$. 
By straightforward computation we have a natural commutative 
diagram
$$
\begin{CD}
0 @>>> T_{[C, P_0]}^*M @>{\varpi^*}>> ((T^*L^\times)/\bC^\times)_{[C, P_0]} 
@>{Z}>> \bC @>>> 0\\
@. @| @VVV @| @.\\
0 @>>> H^0(C; 2K+P_0) @>>> H^0(C; 2K+2P_0)
@>{2\pi\sqrt{-1}\Res_{P_0}}>> \bC @>>> 0
\end{CD}
$$
Here $\Res_{P_0}: H^0(C; 2K+2P_0) \to \bC$ is the residue map of 
quadratic differentials at $P_0$ defined by 
$$
\Res_{P_0}(q_{-2}z^{-2} +  q_{-1}z^{-1} +  q_0 + q_{1}z^1 + 
\cdots)dz^{\otimes 2} = q_{-2},
$$
where $z$ is a complex coordinate centered at $P_0$. 
It is easy to check $q_{-2}$ does not depend on the choice of 
the coordinate $z$. Consequently any  $C^\infty$ family 
$q = \{q(C, P_0)\}_{[C, P_0] \in \Cg}$, $q(C, P_0) \in H^0(C; 2K+2P_0)$
of quadratic differentials parametrized by the space $\Cg$  satisfying
the condition $\Res_{P_0}q(C, P_0) = \frac{1}{2\pi\sqrt{-1}}$ for any
$[C, P_0] \in \Cg$  corresponds to an Ehresmann connection on the
relative tangent bundle $T_{\mathbb{C}_g/\mathbb{M}_g}$. The $(1, 1)$
form $\frac{\sqrt{-1}}{2\pi}\dc q$ on the space $\Cg$ represents the
first Chern class of the bundle $T_{\mathbb{C}_g/\mathbb{M}_g}$
\begin{equation}
\frac{\sqrt{-1}}{2\pi}[\dc q] 
= c_1(T_{\mathbb{C}_g/\mathbb{M}_g})
\in H^2(\Cg; \bR)\label{residue}
\end{equation}
\cite{KawF}.

\section{The Weil-Petersson K\"ahler form}\label{s-wp}

As was shown in \S\ref{s-cotang} the cotangent space 
of the moduli space $\moduli$ at $[C]$ is naturally isomorphic 
to the space of holomorphic quadratic differentials, $H^0(C; 2K)$.
Let $\dvol$ denote the hyperbolic volume form on the Riemann surface $C$. 
It is regarded as a Hermitian metric on the 
relative tangent bundle $T_{\mathbb{C}_g/\mathbb{M}_g}$. 
For any two differentials $q_1, q_2 \in H^0(C; 2K)$ 
the Weil-Petersson pairing
$\langle q_1, q_2\rangle_{\WP}$ is defined by the
integral
\begin{equation*}
\langle q_1, q_2\rangle_{\WP} = \int_Cq_1\overline{q_2}/\!\!\dvol.
\end{equation*}
Here $q_1\overline{q_2}/\!\!\dvol$ is regarded as a $(1, 1)$-form on $C$. 
The pairing induces a Hermitian metric on the moduli space $\Mg$, 
the Weil-Petersson metric. Ahlfors
\cite{Ah1} proved it is K\"ahler.  See \cite{G1} for an alternative
gauge-theoretic proof.  Let $\omega_{\WP}$ denote the K\"ahler form of
the Weil-Petersson metric.
\par
Now recall the original definition of the $i$-th Morita-Mumford classs
$e_i = (-1)^{i+1}\kappa_i$, $i \geq 1$ \cite{Mu} \cite{Mo1}. It is defined to be the integral along the fiber of 
the $(i+1)$-st power of the first Chern class of the relative tangent bundle $T_{\mathbb{C}_g/\mathbb{M}_g}$
\begin{equation}
e_i = (-1)^{i+1}\kappa_i = \int_{\fiber}c_1(T_{\mathbb{C}_g/\mathbb{M}_g})^{i+1} \in H^{2i}(\moduli).
\end{equation}
It is one of the most orthodox ways to obtain 
differrential forms representing the Morita-Mumford classes
to take the integral of powers of the 
hyperbolic Chern form of the relative tangent bundle 
$T_{\mathbb{C}_g/\mathbb{M}_g}$ along the fiber. 
This was carried out by Wolpert \cite{Wo}.
He computed the Chern form $c^{\hyp}_1(T_{\mathbb{C}_g/\mathbb{M}_g})$ of the hyperbolic metric explicitly, 
and he proved
\begin{equation}
\int_{\fiber}{c^{\hyp}_1(T_{\mathbb{C}_g/\mathbb{M}_g})}^2 = \frac1{2\pi^2} \omega_{\WP}
\end{equation}
as differential forms on the moduli space $\moduli$. As a corollary we have
\begin{equation*}
\frac1{2\pi^2} [\omega_{\WP}] = e_1 \in H^2(\moduli; \bR).
\end{equation*}
Furthermore Wolpert \cite{Wo1} gave a description of the Weil-Petersson K\"ahler form in terms of the Fenchel-Nielsen coordinates $(\tau_j, \ell_j)$, $1 \leq j \leq 3g-3$, for any 
pants decomposition of the surface
\begin{equation}
\omega_{\WP} = \sum d\ell_i\wedge d\tau_i.
\end{equation}
Here $\ell_j$ denotes the geodesic length of each simple closed curve in the decomposition, and $\tau_j \in \bR$ the hyperbolic displacement parameter. 
Penner \cite{Pe} described explicitly 
the pullback of $\omega_{\WP}$ to the decorated Teichm\"uller space.  
Goldman \cite{G2} generalized the Weil-Petersson geometry to 
the space of surface group representations in a reductive Lie group.
\par
Now we consider the Lie algebra ${\bold d}$ of complex analytic vector 
fields on the punctured disk $\{z \in \bC; \,\, 0 < \vert z\vert < \epsilon\}$, $0 < \epsilon \ll 1$. The $2$-cochain $\mbox{vir}$ 
on ${\bold d}$ defined by 
\begin{eqnarray}
\mbox{vir}\left(f_1(z)\frac{d}{dz}, f_2(z)\frac{d}{dz}\right) &:= & 
\frac1{2\pi\sqrt{-1}}\oint_{\vert z\vert = 1}\det\begin{pmatrix}
f'_1(z) & f'_2(z) \\ f''_1(z) & f''_2(z)\end{pmatrix} dz\nonumber\\
&=& \frac{\sqrt{-1}}{2\pi}\oint_{\vert z\vert = 1}\det\begin{pmatrix}
f_1(z) & f_2(z) \\ f'''_1(z) & f'''_2(z)\end{pmatrix} dz\nonumber
\end{eqnarray}
is a cocycle and it is called the Virasoro cocycle. 
Its cohomology class generates the second Lie algebra cohomology group
$H^2({\bold d}) = \bC$. \par
Arbarello, De Concini, Kac and Procesi \cite{ADKP} established an 
isomorphism of $H^2({\bold d})$ onto the second cohomology group of 
$\moduli$
\begin{equation}
\nu: H^2({\bold d}) \mathop{\to}^{\cong} H^2(\moduli; \bC)
\end{equation}
induced by the Krichever construction.\par
For a local coordinate $z$ on a Riemann surface one can define 
a local differential operator, or a local complex analytic 
Gel'fand-Fuks $1$-cocycle with values in quadratic differentials by
\begin{equation*}
\nabla^{d/dz}_2: f(z)\frac{d}{dz} \mapsto \frac{1}{6} f'''(z)(dz)^{\otimes 2}
\end{equation*}
\cite[p.666]{Kaw0}. The cocycle $\nabla^{d/dz}_2$ is equivalent to 
a projective structure. In fact, if $w$ is another coordinate, then
\begin{equation*}
\nabla^{d/dw}_2X - \nabla^{d/dz}_2X 
= \mathcal{L}_X\left(\{w, z\} (dz)^{\otimes 2}\right)
\end{equation*}
for any local complex analytic vector field $X$. 
Here $\{w, z\}$ denotes the Schwarzian derivative. 
In particular, the hyperbolic structure on a (hyperbolic) 
Riemann surface defines a global operator $\nabla^{\hyp}_2$. \par
The Krichever construction relates the $2$-cocycle $\mbox{vir}$ 
with the operator $\nabla^{\hyp}_2$. 
By straightforward computation
using the Bers embedding we have
\begin{equation}
\dc\nabla^{\hyp}_2 = 8\omega_{\WP}
\end{equation}
as $(1, 1)$-forms on the moduli space $\moduli$. 
This result, the first variation of the hyperbolic 
structure coincides with $\omega_{\WP}$, was first 
proved by Zograf and Takhtajan \cite[p.310]{ZT1}.

\section{The first Chern form on the Siegel upper
halfspace}\label{s-chern}
The Hodge bundle $\Hodge$ is defined to be the holomorphic 
vector bundle on $\moduli$ whose fiber over $[C]$ is 
the space of holomorphic $1$-forms on $C$
$$
\Hodge = \coprod_{[C] \in \moduli} H^0(C; K).
$$
We write simply $c_1$ for the first Chern class of $\Hodge$
$$
c_1 = c_1(\Hodge) \in H^2(\moduli; \bR).
$$
The bundle $\Hodge$ comes from a symplectic equivariant
vector bundle on the Siegel upper halfspace $\Hg$. 
In fact, the space $\Hg$ can be identified with the space 
of almost complex structures $J$ on the real $2g$-dimensional 
symplectic vector space $(\bR^{2g}, \cdot)$ with the conditions
\begin{eqnarray*}
& Jx\cdot Jy = x\cdot y, \quad \forall x, \forall y \in \bR^{2g},\\
& x\cdot Jx > 0, \quad \forall x \in \bR^{2g}\setminus \{0\}.
\end{eqnarray*}
We have a holomorphic vector bundle $E'_{\Hg}$ on $\Hg$ 
whose fiber over $J$ is the $-\sqrt{-1}$-eigenspace of $J$. 
We have a natural isomorphism of vector bundles
\begin{equation}
T^*\Hg = \Sym^2E'_{\Hg}.
\end{equation}
For each Riemann surface $C$ the Hodge $*$-operator on the 
$1$-forms induces such an almost complex structure on the 
space of real harmonic $1$-forms.
This induces a holomorphic map $\Jac: \moduli \to \Hg/Sp_{2g}(\bZ)$ 
known as the period map in the classical context. 
The pullback of $E'_{\Hg}$ by the map $\Jac$ is exactly 
the Hodge bundle $\Hodge$. \par
Thus the cohomology class $c_1$ can be regarded as an integral cohomology class of the Siegel modular group $Sp_{2g}(\bZ)$,
$c_1 \in H^2(Sp_{2g}(\bZ); \bZ)$. 
Meyer \cite{Me1} proved that the cohomology class of 
the Meyer cocycle 
is equal to $4c_1 \in H^2(Sp_{2g}(\bZ); \bZ)$. From the 
Grothendieck-Riemann-Roch formula, or equivalently 
the Atiyah-Singer index theorem for families, 
it follows that
\begin{equation}
\frac1{12}e_1 = c_1 \in H^2(\moduli; \bR).
\label{1/12}\end{equation}
To describe a canonical $2$-form representing $c_1(E'_{\Hg})$ 
we consider the quotient vector bundle $E''_{\Hg} := (\Hg\times\bC^{2g})/E'_{\Hg}$, and the family of projections 
$\pi = \{\pi_J\}_{J \in \Hg}$ on $\bC^{2g}$, $\pi_J := \frac{1}{2}(1 - \sqrt{-1}J)$, parametrized by $\Hg$. 
Then $\{\pi_J \circ d\}_{J \in \Hg}$ is a covariant derivative 
$\nabla$ of type $(1, 0)$ on the bundle $E''_{\Hg} \cong \coprod_{J \in \Hg}\operatorname{Image}\pi_J$, 
whose curvature form $R^\nabla$ is given by 
\begin{equation}
R^\nabla = \pi(\dd\pi)(\dc\pi).
\end{equation}
The $2$-form $c_1(\nabla)$ defined by $c_1(\nabla) = \frac{\sqrt{-1}}{2\pi}\trace R^\nabla$ represents $c_1(E'_{\Hg})$. 
Let $J_\alpha(t) \in \Hg$, $\vert t\vert \ll 1$, $\alpha = 1, 2$, 
be $C^\infty$ paths on $\Hg$ with $J_1(0) = J_2(0) = J$. 
Then, one can compute 
\begin{equation}
c_1(\nabla)_J = \frac1{8\pi}
\trace(\mathop{J_1}^{\centerdot} J\mathop{J_2}^\centerdot).
\label{chern}
\end{equation}
In the next section we prove Rauch's variational formula 
to obtain the pullback of $c_1(\nabla)_J$ by the period map 
$\Jac$ explicitly. 

\section{Rauch's variational formula}\label{s-rau}

Rauch's variational formula describes the differential of 
the period map $\Jac$. Let $C$ be a compact Riemann surface 
of genus $g$. We denote by $H$ the real first homology group 
$H_1(C; \bR)$. Consider the map $H^* = H^1(C; \bR) \to \Omega^1(C)$ 
assigning to each cohomology class the harmonic $1$-form 
representing it. The map can be regarded as an $H$-valued $1$-form 
$\omegaone \in \Omega^1(C)\otimes H$. \par
Let $\{X_i, X_{g+i}\}^g_{i=1}$ be a symplectic basis of 
$H_\bC = H_1(C; \bC)$
$$
X_{i}\cdot X_{g+j} = \delta_{ij}, \quad
X_{i}\cdot X_{j} = X_{g+i}\cdot X_{g+j} = 0, \quad
1 \leq i, j \leq g,
$$
and  $\{\xi_i, \xi_{g+i}\}^g_{i=1} \subset \Omega^1(C)$ 
the basis of the harmonic $1$-forms dual to  $\{X_i, X_{g+i}\}^g_{i=1}$. 
Then we have 
\begin{equation*}
\omegaone = \sum^g_{i=1}\xi_{i}X_{i} + \xi_{g+i}X_{g+i} 
\in \Omega^1(C)\otimes H_\bC.
\end{equation*}
In particular, if $\{\psi_i\}^g_{i=1} \subset H^0(C; K)$ is an orthonormal basis
\begin{equation}
\frac{\sqrt{-1}}{2}\int_C\psi_i\wedge\overline{\psi_j} 
= \delta_{ij}, \quad 1 \leq i, j \leq g, 
\label{onb}
\end{equation}
then we obtain
\begin{equation}
\omegaone = \sum^g_{i=1} \psi_iY_i + \overline{\psi_i}\overline{Y_i},\label{psiY}
\end{equation}
where $\{Y_i, Y_{g+i}\}^g_{i=1} \subset H_\bC$ is the dual basis of the symplectic basis $\{[\psi_i], \frac{\sqrt{-1}}{2}[\overline{\psi_i}]\}^g_{i=1}$ of $H^*_\bC = H^1(C; \bC)$. Since the complete linear system of the canonical divisor on the complex algebraic curve $C$ has no basepoint, the $2$-form 
\begin{equation}
B = \frac{1}{2g}\omegaone\cdot\omegaone = 
\frac{\sqrt{-1}}{2g}\sum^g_{i=1}\psi_i\wedge\overline{\psi_i} 
\in \Omega^2(C)\label{B}
\end{equation}
is a volume form on $C$. 
\par
Now we recall the Hodge decomposition of the $1$-forms on $C$. We have an exact sequence
\begin{equation*}
0 \to \bC \to \Omega^0(C) \mathop{\to}^{d*d} \Omega^2(C)
\mathop{\to}^{\int_C} \bC \to 0.
\end{equation*}
The vector space $\bC$ on the left side means the constant
functions. 
A Green operator $\Psi: \Omega^2(C) \to \Omega^0(C)$ is a 
linear map satisfying the property 
$$
d* d\Psi \Omega = \Omega 
$$
for any $\Omega \in \Omega^2(C)$ with $\int_C\Omega = 0$.
In this article we use two sorts of Green operators $\hatPhi = 
\hatPhi_C$ and $\Phi = \Phi^{(C, P_0)}$. The former is characterized by the conditions
\begin{equation}
d* d\hatPhi(\Omega) = \Omega - (\int_C\Omega)B \quad \mbox{and}\quad  \int_C\hatPhi(\Omega)B = 0
\label{hatPhi}\end{equation}
for any $\Omega \in \Omega^2(C)$. 
Let $\delta_{P_0}: C^\infty(C) \to \bC$,
$f \mapsto f(P_0)$, be the delta current on $C$ at 
the point $P_0$. 
We define the latter $\Phi$ to be 
a linear map with values in $\Omega^0(C)/\bC$ instead of
$\Omega^0(C)$. Then the operator $d\Phi: \Omega^2(C) \to \Omega^1(C)$ makes sense, and the operator $\Phi$ 
is defined by the condition 
$$
d* d\Phi \Omega = \Omega - \left(\int_C\Omega\right)\delta_{P_0}
$$
for any $\Omega \in \Omega^2(C)$.
\par
Any Green operator $\Psi$ induces the Hodge
decomposition of the $1$-currents 
\begin{equation}
\varphi = \harmonic\varphi + d\Psi d\ast \varphi + \ast d\Psi d \varphi\label{decomp}
\end{equation}
for any $\varphi \in \Omega^1(C)$, where $\harmonic: \Omega^1(C) \to \Omega^1(C)$ is the
harmonic projection on the $1$-currents on $C$.
\par
In the setting of \S\ref{s-cotang} the first variation of 
$\omegaone$ is given by
\begin{equation}
\omegaonedot = -d\Psi d* S\omegaone.
\label{omegaonedot}
\end{equation}
In fact, differentiating $d*\omegaone = 0$, we get
$$
d*\omegaonedot = -d\mathop{*}^{\centerdot}\omegaone = -d*S\omegaone.
$$ 
Since ${f^t}^*\omegaone$ is cohomologous to $\omegaone$, 
we have some function $u$ such that $\omegaonedot = du$. Hence from (\ref{decomp}) we obtain
$$
\omegaonedot = d\Psi d*\mathop{\omegaone}^{\centerdot} = -d\Psi d*\omegaone,
$$
as was to be shown.

\begin{theorem}[Rauch] The diagram
$$
\begin{CD}
T^*_{[C]}\moduli @<{(d\operatorname{Jac})^*}<< T^*_{[\operatorname{Jac}(C)]}\Hg/Sp_{2g}(\bZ)\\
@| @|\\  
H^0(C; 2K) @<{2\sqrt{-1}\text{(multiplication)}}<< \operatorname{Sym}^2H^0(C; K).
\end{CD}
$$
commutes. 
Here the lower horizontal arrow maps $\psi_1\otimes\psi_2$ 
to the quadratic differential $2\sqrt{-1}\psi_1\psi_2$ for 
any $1$-forms $\psi_1$ and $\psi_2 \in H^0(C; K)$. 
\end{theorem}
\begin{proof}
The integral $\int_C*\omegaone\wedge\omegaone \in H\otimes H = H^* \otimes H = \Hom(H, H)$ coincides with 
the almost complex structure on $H = H_1(C; \bR)$ 
induced by the Hodge $*$-operator. 
Since $\omegaone$ is harmonic and $\omegaonedot$ is $d$-exact by (\ref{omegaonedot}), we have 
$$
\int_C*\omegaone\wedge\omegaonedot = 
- \int_C\omegaone\wedge*\omegaonedot = 0.
$$
Hence
\begin{eqnarray*}
& \left(\int_C*\omegaone\wedge\omegaone\right)^\centerdot
= \int_C\stardot\omegaone\wedge\omegaone = 
\int_C(*S\omegaone)\wedge\omegaone\\
&=  2\sqrt{-1}\int_C\omegaone'\omegaone'\mudot
- 2\sqrt{-1}\overline{\left(\int_C\omegaone'\omegaone'\mudot\right)}.
\end{eqnarray*}
This proves the theorem.
\end{proof}
Substituting the theorem into the formula (\ref{chern}) 
we have 
\begin{corollary}
$$
\Jac^*c_1(\nabla) = \frac1{8\pi\sqrt{-1}}\sum^g_{i, j= 1}
\psi_i\psi_j\otimes\overline{\psi_i}\overline{\psi_j} 
\in T^*_{[C]}\moduli\otimes \overline{T^*_{[C]}\moduli}.
$$
Here $\{\psi_i\}^g_{i=1} \subset H^0(C; K)$ is any orthonormal 
basis (\ref{onb}).
\end{corollary}

The elementary polynomials $\sigma_1, \dots, \sigma_g$ in indeterminates 
$x_1, \dots, x_g$ are 
given by $\prod^g_{i=1}(t-x_i) = t^g + \sum^g_{k=1}(-1)^k\sigma_k t^{g-k}$. The equation $\sum^g_{i=1} {x_i}^m = s_m(\sigma_1, \dots, \sigma_g)$ defines the $m$-th Newton polynomial $s_m$. 
The $m$-th Newton class of the Hodge bundle 
$\Lambda = \Lambda_{\moduli}$ is defined by
$$
s_m(\Lambda) = s_m(c_1(\Lambda), \dots, c_g(\Lambda)) \in H^{2m}(\moduli; \bR),
$$
where $c_k(\Lambda)$ is the $k$-th Chern class of 
the bundle $\Lambda$. \par
The complex conjugate $\overline{\Lambda}$ satisfies 
$s_m(\overline{\Lambda}) = (-1)^ms_m(\Lambda)$. 
Since $\Lambda\oplus\overline{\Lambda}$ is a flat vector bundle on $\moduli$ whose fiber over $[C]$ is the homology group $H_1(C; \bC)$, we have
$$
s_{2n}(\Lambda) = \frac12s_{2n}(\Lambda\oplus\overline{\Lambda}) = 0.
$$
From the Grothendieck-Riemann-Roch formula or equivalently 
the Atiyah-Singer index theorem for families, 
it follows that 
\begin{equation}
e_{2n-1} = (-1)^{n-1}\frac{2n}{B_{2n}}s_{2n-1}(\Lambda) 
\in H^{4n-2}(\moduli; \bR).
\end{equation}
Here $B_{2n}$ is the $n$-th Bernoulli number.
In the case $n=1$ it is exactly the formula (\ref{1/12}).  \par
Hence the Hodge bundle yields all the
odd Morita-Mumford classes,
 but not the even ones.
To get all the Morita-Mumford classes we 
introduce a higher analogue of the period map,
as will be discussed in the succeeding sections.

\section{The Earle class and the twisted Morita-Mumford classes}\label{s-earle}

Let $\Sigma_g$ be a closed oriented $C^\infty$ surface of genus $g$, $p_0 \in \Sigma_g$ a point, and $v_0 \in T_{p_0}\Sigma_g\setminus\{0\}$ a non-zero tangent vector at the point $p_0$. 
We denote by $\Mg$, $\Mgstar$ and $\Mgone$ the mapping 
class groups for the surface $\Sigma_g$, the pointed 
surface $(\Sigma_g, p_0)$ and the triple $(\Sigma_g, p_0, v_0)$ respectively. 
They are the orbifold fundamental groups of the spaces $\moduli$, $\Cg$ and $\TCM$. 
The fundamental group $\pi_1(\Sigma_g, p_0)$ is naturally embedded into the group $\Mgstar$ \cite{MoH}. \par
By abuse of notation let $H$ denote the real first homology group of $\Sigma_g$, $H_1(\Sigma_g; \bR)$, on which 
the mapping class groups act in an obvious way. 
The module $H$  can be interpreted as a flat vector bundle on the moduli 
space $\moduli$. In 1978 Earle \cite{E} constructed an explicit $1$-cocycle $\psi: \Mgstar \to H$ such that 
$(2-2g)\psi$ has values in $H_1(\Sigma_g; \bZ)$, 
and $\psi\vert_{\pi_1(\Sigma_g)}$ is equal to the abelianization map of the group $\pi_1(\Sigma_g)$.
Later Morita \cite{Mo3} independently discovered 
a cohomology class $k \in H^1(\Mgstar; H_1(\Sigma_g; \bZ))$ which is equal to $[(2-2g)\psi]$. Furthermore he proved 
\begin{equation}
H^1(\Mgstar; H_1(\Sigma_g; \bZ)) = \bZ k \cong \bZ
\end{equation}
for $g \geq 2$. The author would like to propose 
the class $k$ should be called {\it the Earle class}. \par
The square of the class $k$ is related to the first Morita-Mumford class 
$e_1 = \kappa_1$ through the intersection pairing
\begin{equation}
m: H \otimes H = H_1(\Sigma_g; \bR) \otimes H_1(\Sigma_g; \bR) \to \bR.
\label{intersect}
\end{equation}
Morita \cite{Mo4} proved 
\begin{equation}
m_*(k^{\otimes 2}) = -e_1 + 2g(2-2g)e \in H^2(\Mgstar).
\label{contr}\end{equation}
Here $e$ is the first Chern class of the relative tangent bundle $c_1(T_{\mathbb{C}_g/\mathbb{M}_g}) \in H^2(\Cg) = H^2(\Mgstar)$. \par
These phenomena have a higher analogue. 
The twisted Morita-Mumford class  $m_{i, j} \in H^{2i+j-2}(\Mgone; \Lambda^jH)$, $i, j \geq 0$,  was introduced in 
\cite{Kaw1}. We have $m_{1,1} = k$ and $m_{i+1, 0} = e_i$, $i\geq 1$. All the cohomology classes on the mapping class groups with trivial coefficients (even in the unstable range) obtained from any products of the twisted Morita-Mumford classes by contracting the coefficients using the intersection pairing are exactly the polynomials in the Morita-Mumford classes \cite{KM}. \par
This fact is closely related to the Johnson homomorphisms
on the mapping class group. 
The fundamental group $\pi_1(\Sigma_g, p_0, v_0) 
= \pi_1(\Sigma_g\setminus\{p_0\}, v_0)$ with 
tangential basepoint $v_0$ is a free group of rank 
$2g$. Let $\Gamma_k$, $k \geq0$, denote the lower 
central series of the free group $\pi_1(\Sigma_g, p_0, v_0)$. 
We have $\Gamma_0 = \pi_1(\Sigma_g, p_0, v_0)$ 
and $\Gamma_{k+1} = [\Gamma_k, \Gamma_0]$ for $k \geq 0$. 
The quotient $\Gamma_1/\Gamma_2$ is naturally isomorphic to
$\bigwedge^2H_1(\Sigma_g; \bZ) \subset  \bigwedge^2H$. Let
$\mathcal{I}_{g, 1}$ be the Torelli group, that is, 
the kernel of the natural action of 
$\mathcal{M}_{g, 1}$ on the homology group $H_1(\Sigma_g; \bZ)$. 
For any $\varphi \in \mathcal{I}_{g, 1}$ and $\gamma \in \Gamma_0$, the
difference $\gamma^{-1}\varphi(\gamma)$ belongs to $\Gamma_1$  from the
definition of $\mathcal{I}_{g, 1}$. Hence we can define a homomorphism
$$
\tau_1(\varphi): H_1(\Sigma_g; \bZ) \to {\bigwedge}^2H_1(\Sigma_g; \bZ), \quad [\gamma]\mapsto 
\gamma^{-1}\varphi(\gamma) \bmod \Gamma_2.
$$
It is easy to check this induces a homomorphism 
$\tau_1: \mathcal{I}_{g, 1} \to H^*\otimes {\bigwedge}^2H \cong 
H\otimes {\bigwedge}^2H$. The last isomorphism comes from Poincar\'e
duality.  Johnson \cite{J} proved the image $\tau_1(\mathcal{I}_{g, 1})$
is included in ${\bigwedge}^3H$. The homomorphism $\tau_1$ is called  the
first Johnson homomorphism.  Morita \cite{Mo6} proved there exists a
unique  cohomology class $\tilde k \in H^1(\Mgone; {\bigwedge}^3H)$ which
restricts to $\tau_1$ on the Torelli group 
$\mathcal{I}_{g, 1}$. We call it 
the extended first Johnson homomorphism. 
See \cite[\S7]{MoH} for
more information on the Johnson homomorphisms. \par
The class $\frac1{6}m_{0, 3}$ is equal to the extended first 
Johnson homomorphism $\tilde k: \Mgone \to \bigwedge^3H$ \cite{KM}. 
Each of the Morita-Mumford classes is obtained from some power of $\tilde
k$ by contracting the coefficients using the intersection pairing $m$ 
\cite{MoF}. 
Conversely for any $\Symp$-module $V$ and any $\Symp$-homomorphism $f: (\bigwedge^3H)^{\otimes m} \to V$ induced by the intersection pairing, the cohomology class $f_*({\tilde k}^{\otimes m})$ is a polynomial in the twisted Morita-Mumford class \cite{KM}. An extension of the second Johnson homomorphism to the whole mapping class group provides a fundamental relation among the twisted Morita-Mumford classes \cite{Kaw2}. 
In the next section we introduce a flat connection on a vector bundle on the space $\TCM$, whose holonomy is an extension of the Johnson homomorphisms to the whole mapping class group $\Mgone$. 

\section{A higher analgue of the period map}\label{s-hme}

A complex-analytic counterpart of the first Johnson homomorphism is the (pointed) harmonic volume introduced by Harris \cite{BH} \cite{Pu}. It is a real analytic section of a fiber bundle on the moduli $\Cg$ whose fiber over $[C, P_0]$ is $(\bigwedge^3H_1(C; \bZ))\otimes (\bR/\bZ)$. 
The first variation of the (pointed) harmonic volumes is a twisted $1$-form representing the cohomology class $[\tilde k]$ \cite {Kaw3}. \par
To obtain ``canonical" differential forms representing all the 
twisted Morita-Mumford classes and their higher relations, we construct a higher analogue of the classical period map
and the harmonic volume, the harmonic Magnus expansion
$\theta: {\mathcal T}_{g, 1} \to \Theta_{2g}$
\cite{Kaw3}. 
The space $\mathcal T_{g, 1} = \widetilde{\TCM}$ is  Teichm\"uller space of 
{\it triples} $(C, P_0, v)$ of genus $g$. 
Here $C$ is a compact Riemann surface 
of genus $g$, $P_0 \in C$, and $v$ a non-zero tangent 
vector of $C$ at $P_0$ as in \S\ref{s-cotang}. 
For any triple $(C, P_0, v)$ one can define the fundamental 
group of the complement $C\setminus \{P_0\}$ with the tangential basepoint $v$ denoted by $\pi_1(C, P_0, v)$, which is a free group of rank $2g$.
The space $\Theta_n$ is the set of all 
Magnus expansions of the free group $F_n$ of rank $n \geq 2$ 
in a wider sense stated as follows.\par
We denote by $H$ 
the first real homology group of the group $F_n$, 
$H_1(F_n; {\mathbb R})$, 
 $H^*$ 
the first real cohomology group of $F_n$, 
$H^1(F_n; {\mathbb R})$,
and $[\gamma] \in H$ 
the homology class of $\gamma \in F_n$. 
The completed tensor algebra generated by $H$, 
$\widehat{T} = \widehat{T}(H) = {\prod}^\infty_{m=0} 
H^{\otimes m}$, 
has a decreasing filtration of two-sided ideals
$\{\widehat{T}_p\}_{p\geq 1}$ defined by 
$\widehat{T}_p = {\prod}_{m\geq p} H^{\otimes m}$.
The subset $1+{\widehat{T}}_1$ is 
a subgroup of the multiplicative group of 
the algebra ${\widehat{T}}$.
We call a map $\theta: F_n \to 1 +
{\widehat{T}}_1$ {\it a Magnus expansion of the free group $F_n$}
 in a wider sense \cite{Kaw2}, 
if $\theta: F_n \to 1 + {\widehat{T}}_1$ is a
group homomorphism, and if 
$\theta(\gamma) \equiv 1 + [\gamma] \pmod{{\widehat{T}}_2}$
for any $\gamma \in F_n$.
One can endow the set of all Magnus expansions
$\Theta_n$ with a natural strucure of a (projective limit 
of) real analytic  manifold(s). 
A certain (projective limit of) Lie group(s)
$\operatorname{IA}(\widehat{T})$ acts  on $\Theta_n$ in a free and 
transitive way.  This induces a series of $1$-forms 
$\eta_p \in \Omega^1(\Theta_n)\otimes H^*\otimes 
H^{\otimes (p+1)}$, $p \geq 1$, the Maurer-Cartan forms 
of the action of $\operatorname{IA}(\widehat{T})$, 
which are invariant under a natural action of the automorphism 
group of the group $F_n$, $\operatorname{Aut}(F_n)$.  
The Maurer-Cartan formula $d\eta = \eta\wedge\eta$ allows us to regard the forms $\eta_p$ as an equivariant flat connection on the vector bundle $\Theta_n\times H^*\otimes \widehat{T}_2$.
The holonomy of the connection is an extension of all the Johnson homomorphisms to the whole group 
$\operatorname{Aut}(F_n)$. 
The $1$-forms $\eta_p$ represent the twisted Morita-Mumford classes on the group $\operatorname{Aut}(F_n)$ \cite{Kaw2} \cite{Kaw3}. 
\par
Let $(C, P_0, v)$ be a triple of genus $g$. 
From now on we denote by $H$ the real first homology group $H_1(C; \bR)$. 
As in \S\ref{s-rau} we denote by $\delta_{P_0}: C^\infty(C) \to {\mathbb R}$, $f\mapsto
f(P_0)$, the delta $2$-current on $C$ at $P_0$. 
Then there exists a $\widehat{T}_1$-valued $1$-current $\omega \in
\Omega^1(C)\otimes
\widehat{T}_1$, satisfying the following 3 conditions
\begin{enumerate}
\item[(1)] $d\omega = \omega\wedge\omega - I\cdot \delta_{P_0}$, where 
$I \in H^{\otimes 2}$ is the intersection form.
\item[(2)] The first term of $\omega$ is equal 
to $\omega_{(1)} \in \Omega^1(C)\otimes H$ introduced in \S\ref{s-rau}.
\item[(3)] $\int_C(\omega - \omega_{(1)})\wedge*\varphi = 0$ for any
closed
$1$-form
$\varphi$ and each $p\geq 2$. 
\end{enumerate}
Using Chen's iterated integrals \cite{Chen}, 
 we can define a Magnus expansion
$$
\theta = \theta^{(C, P_0, v)}: \pi_1(C, P_0, v) 
\to 1 + \widehat{T}_1(H_1(C; {\mathbb R})), \quad
[\ell] \mapsto 1 + \sum^\infty_{m=1} \int_\ell
\overbrace{\omega\omega\cdots\omega}^{m}.
$$

Let a point $p_0 \in \Sigma_g$ and a non-zero tangent vector $v_0 \in T_{p_0}\Sigma_g \setminus\{0\}$ be fixed
as in \S\ref{s-earle}. 
Moreover we fix an isomorphism $\pi_1(\Sigma_g, p_0, v_0) 
\cong F_{2g}$. A marking $\alpha$ of a triple $(C, P_0, v)$ 
is an orientation-preserving diffeomorphism of $\Sigma_g$ 
onto $C$ satisfying the conditions $\alpha(p_0) = P_0$ 
and $(d\alpha)_{p_0}(v_0) = v$. 
For any marked triple $[(C, P_0, v), \alpha]$ we define a 
Magnus expansion of the free group $F_{2g}$ by 
$$
F_{2g}\cong \pi_1(\Sigma_g, p_0, v_0) \mathop\rightarrow^{\alpha_*}
\pi_1(C, P_0, v) \mathop{\longrightarrow}^{\theta^{(C, P_0, v)}}
1 + \widehat{T}_1(H_1(C; {\mathbb R}))
\mathop\rightarrow^{{\alpha_*}^{-1}} 1 + \widehat{T}_1. 
$$
Consequently, the Magnus expansions $\theta^{(C, P_0, v)}$ 
for all the triples $(C, P_0, v)$ 
define a canonical real analytic map 
$\theta: \widetilde{\TCM}={\mathcal T}_{g, 1}\to \Theta_{2g}$, 
which we call {\it the harmonic Magnus expansion on the universal family
of  Riemann surfaces}.  The pullbacks of the Maurer-Cartan forms
$\eta_p$  define a flat connection on a vector bundle on the space
$\TCM$, and give the canonical differential forms representing the
Morita-Mumford classes and their higher relations.

\begin{theorem}[\cite{Kaw3}] For any $[C, P_0, v, \alpha]
\in {\mathcal T}_{g, 1}$ we have 
$$
(\theta^*\eta)_{[C, P_0, v, \alpha]} = 2\Re(N(\omega'\omega') -
2{\omega_{(1)}}'{\omega_{(1)}}') 
\in T^*_{[C, P_0, v, \alpha]}{\mathcal T}_{g, 1}\otimes\widehat{T}_3.
$$
Here $N: \widehat{T}_1 \to \widehat{T}_1$ is defined by $
N\vert_{H^{\otimes m}} =\sum^{m-1}_{k=0}\begin{pmatrix}
1& 2& \cdots & m-1 & m\\
2& 3& \cdots & m & 1
\end{pmatrix}^k
$, and the meromorphic quadratic differential $N(\omega'\omega')$ is 
regarded as a $(1, 0)$-cotangent vector at $[C, P_0, v, \alpha]
\in {\mathcal T}_{g, 1}
$ in a natural way.
\end{theorem}

The third homogeneous term $N(\omega'\omega')_{(3)} = N(\omega'_{(1)}\omega'_{(2)} + \omega'_{(2)}\omega'_{(1)})$ is the first variation of the (pointed) harmonic volumes of pointed Riemann surfaces. 
It represents the extended first Johnson homomorphism $\tilde k$. The higher terms provide higher relations among the twisted Morita-Mumford classes. 
Hence all of the Morita-Mumford classes are represented by some algebraic combinations 
of $N(\omega'\omega')$. \par
The second term coincides with 
$2{\omega_{(1)}}'{\omega_{(1)}}'$, which is exactly the first variation of
the  period matrices given by Rauch's formula in \S\ref{s-rau}. Hence we may regard the harmonic Magnus expansion as a higher analogue of the classical period map $\Jac$. \par

\section{Secondary objects on the moduli space}\label{s-fiber}

The determinant of the Laplacian acting on the space of $k$-differentials on Riemann surfaces is a 
`secondary' object on the moduli space. Zograf and Takhtajan \cite{ZT2} proved that it yields the difference on the moduli space of compact Riemann surfaces, $\moduli$, between a multiple of the Weil-Petersson form $\omega_{\WP}$ and the Chern form of the Hodge line bundle for the $k$-differentials induced by the hyperbolic metric. Moreover, they studied analogous phenomena for punctured Riemann surfaces to introduce  their K\"ahler metric, the Zograf-Takhtajan metric, 
on the moduli space of punctured Riemann surfaces \cite{ZT3}. \par
In this section we discuss other secondary objects, 
which come from the higher analogue of the period map 
introduced in \S\ref{s-hme}.
Now we can obtain explicit $2$-forms from the connection form $N(\omega'\omega')$ on $\TCM$, $e^J$ on $\Cg$ and 
$e^J_1$ on $\moduli$. Consider the quadratic differential 
$\eta'_2$ defined by 
$$
\eta'_2 = N(\omega'\omega')_{(4)} \in H^0(C; 2K+2P_0)\otimes H^{\otimes 4},
$$
which satisfies 
$$
\frac{1}{2g(2g+1)}\Res_{P_0}\left((m\otimes m)(\eta'_2)\right) = -\frac1{8\pi^2}.
$$
Here $m$ is the intersection pairing $m: H\otimes H \to \bR$ as in (\ref{intersect}). 
We define 
$$
e^J = \frac{-2}{2g(2g+1)}\dc((m\otimes m)(\eta'_2))
\in \Omega^{1, 1}(\Cg).
$$
From (\ref{residue}) $e^J$ represents the first Chern class of the relative tangent bundle
$$
[e^J] = e = c_1(T_{\mathbb{C}_g/\mathbb{M}_g}) 
\in H^2(\Cg; \bR). 
$$
We obtain a twisted $1$-form $\eta^H_1 \in \Omega^1(\Cg; H)$
representing the Earle class $k$ by contracting the 
coefficients of $\eta'_1 = N(\omega'\omega')_{(3)}$. 
By (\ref{contr}) $m(\eta^H_1)^{\otimes 2} \in \Omega^{1, 1}(\Cg)$ represents $-e_1+2g(2-2g)e$. So we define
$$
e^J_1 = -m(\eta^H_1)^{\otimes 2} + 2g(2-2g)e^J
$$
which can be regarded as a $(1, 1)$-form on $\moduli$ \cite[\S8]{Kaw3}. \par
Hain and Reed \cite{HR} already constructed the same form $e^J_1$ in a Hodge-theoretical context. They applied the following lemma to $\frac{1}{12}e^J_1 - \Jac^*c_1(\nabla)$ to get a 
function $\beta_g \in C^\infty(\moduli; \bR)/\bR$, 
the Hain-Reed function, a secondary object on the moduli space $\moduli$.
\begin{lemma}\label{d-d}
Let $M$ be a connected complex orbifold with $H^0(M; \mathcal{O}) = \bC$ and $H^1(M; \bC) = H^1(M; \mathcal{O}) = 0$. If a real $C^{\infty}$ $(1, 1)$-form $\psi$ is $d$-exact, 
then there exists a real-valued function $f \in C^\infty(M; \bR)$ such that $
\psi = \frac{\sqrt{-1}}{2\pi}\dd\dc f$. 
Such a function $f$ is unique up to a constant.
\end{lemma}
Here we remark all the holomorphic functions on $\moduli$ are constants provided $g \geq 3$. In fact, each of the boundary component of the Satake compactification of $\moduli$ is of complex codimension $\geq 2$. 
The vanishing of the first cohomology follows from (\ref{b1}). See \cite{Mu1}.
Hain and Reed also studied the asymptotic behavior of the  function $\beta_g$ towards the boundary of the Deligne-Mumford compactification $\overline{\moduli}^{\text{DM}}$
\cite{HR}. \par
We have another `secondary' phenomenon arround the $2$-forms $e^J$ 
and $e_1^J$ \cite{Kaw4}.
Let $B = \frac1{2g}\omegaone\cdot\omegaone$ be the volume form in
(\ref{B}).   On any pointed Riemann surface $(C, P_0)$ there exists a
function $h = h_{P_0} = -\hatPhi(\delta_{P_0})$ with $d\ast dh = B -
\delta_{P_0}$ and $\int_ChB = 0$. The function $G(P_0, P_1)  :=
\exp(-4\pi h_{P_0}(P_1))$ is just the Arakelov-Green function. We regard
$G$ a function on the fiber product $\Cg\times_\moduli\Cg$  and define
the $(1, 1)$-form $e^A$ on $\Cg$ by 
$$
e^A := \frac1{2\pi\sqrt{-1}}\dd\dc\log G\vert_{\text{diagonal}} \in \Omega^{1, 1}(\Cg)
$$
representing the Chern class $e = c_1(T_{\mathbb{C}_g/\mathbb{M}_g})$. In fact, the normal bundle of the diagonal map $\Cg \to \Cg\times_\moduli\Cg$ is exactly the relative tangent bundle $T_{\mathbb{C}_g/\mathbb{M}_g}$. 
\par
Furthermore we introduce an explicit real-valued function 
$a_g$ on $\moduli$ by
\begin{eqnarray}
a_g(C) &:=& \int_C\omegaone\cdot\hatPhi(\omegaone\wedge\omegaone)
\cdot\omegaone\label{eqn:a0},
\end{eqnarray}
where $\hatPhi$ is the Green operator introduced in (\ref{hatPhi}). 
By (\ref{psiY}) we have 
\begin{equation}
a_g(C) = - \sum^g_{i,j=1} \int_C\psi_i\wedge\overline{\psi_j}
\hatPhi(\overline{\psi_i}\wedge\psi_j).
\end{equation}
We have $a_g(C) > 0$ if $g \geq 2$. 
Then comparing $\dd a_g$ with $\eta'_2$ as explicit quadratic differentials, we obtain 
\begin{equation}
e^A - e^J = \frac{-2\sqrt{-1}}{2g(2g+1)}\dd\dc a_g.
\end{equation}
On the other hand, the integral along the fiber
$$
e^F_1 := \int_{\text{fiber}} (e^J)^2 \in \Omega^{1, 1}(\moduli)
$$
also represents the first Morita-Mumford class $e_1$. 
By straightforward computation on $\dd\dc a_g$ we deduce
\begin{theorem}[\cite{Kaw4}]
$$
e^A - e^J =  \frac{-2\sqrt{-1}}{2g(2g+1)}\dd\dc a_g = \frac1{(2-2g)^2}({e_1}^F -
{e_1}^J). 
$$
\end{theorem}
The function $a_g(C)$ is also a secondary object on the moduli space $\moduli$, and it defines 
a conformal invariant of the compact Riemann surface $C$, 
but the author does not know any of its further properties.

\end{document}